\documentclass[a4paper]{article}
\usepackage{amsmath}
\usepackage{amssymb}
\usepackage{amsfonts}
\usepackage{theorem}
\usepackage{blkarray}
\usepackage{arydshln}
\usepackage{multirow}
\usepackage{graphicx}

\newcommand{\bigzerou}{%
\smash{\lower1.7ex\hbox{\bg 0}}}
\newtheorem{theorem}{Theorem}

\newtheorem{defi}{Definition}

\newtheorem{lem}{Lemma}

\newcommand{\ba}{\begin{eqnarray}}
\newcommand{\ea}{\end{eqnarray}}
\newcommand{\ban}{\begin{eqnarray*}}
\newcommand{\ean}{\end{eqnarray*}}
\newcommand{\no}{\nonumber}

\def\d{{\partial}}
\newcommand{\mapright}[1]{%
\smash{\mathop{%
\hbox to 1.0cm{\rightarrowfill}}\limits^{#1}}}
\newcommand{\mapleft}[1]{%
\smash{\mathop{%
\hbox to 1.3cm{\leftarrowfill}}\limits^{#1}}}
\allowdisplaybreaks[4]
\begin{document}
\title{
\begin{flushright}
  \begin{minipage}[b]{5em}
    \normalsize
    ${}$      \\
  \end{minipage}
\end{flushright}
{\bf Generalized Hypergeometric Functions for Degree $k$ Hypersurface in $CP^{N-1}$ and Intersection Numbers of Moduli Space of Quasimaps from $CP^{1}$ with Two Marked Points to $CP^{N-1}$ }}
\author{Masao Jinzenji${}^{(1)}$, Kohki Matsuzaka${}^{(2)}$ \\
\\ 
${}^{(1)}$ \it Department of Mathematics,  \\
\it Okayama University \\
\it  Okayama, 700-8530, Japan\\
\\
${}^{(2)}$\it Faculty of Integrated Media, \\
\it Ikueikan University \\
\it  Wakkanai, 097-0013, Japan\\
\\
\it e-mail address: 
\it\hspace{0.0cm}${}^{(1)}$ pcj70e4e@okayama-u.ac.jp \\
\it\hspace{2.7cm}${}^{(2)}$ kohki@ikueikan.ac.jp } 
\maketitle
\begin{abstract}
In this paper, we derive the generalized hypergeometric functions used in mirror computation of degree $k$ hypersurface in $CP^{N-1}$ as generating functions of 
intersection numbers of the moduli space of quasimaps from $CP^{1}$ with two marked points to $CP^{N-1}$.   
\end{abstract}
\section{Introduction}
.In this paper, we discuss the following two (intersection) numbers defined as values of residue integrals.
\begin{defi}
\ba
&&w(\sigma_{(N-k)d+j-1}({\cal O}_{h^{N-2-j}}){\cal O}_{h^{0}})_{0,2}:=\no\\
&&\frac{1}{(2 \pi \sqrt{-1})^{d+1}} \oint_{C_{0}} \frac{dz_{0}}{(z_{0})^{N}} \oint_{C_{1}} \frac{dz_{1}}{(z_{1})^{N}} \dots \oint_{C_{d}} \frac{dz_{d}}{(z_{d})^{N}}(z_{0})^{N-2-j} (z_{1} - z_{0})^{(N-k)d+j-1} \biggl(\prod_{l=1}^{d} e^{k}(z_{l-1},z_{l}) \biggr)\no\\
&&\times\prod_{l=1}^{d-1} \frac{1}{kz_{l} (2z_{l} - z_{l-1} - z_{l+1})} \hspace{7cm}(N>k\geq 1).
\label{residue1}
\ea
\ba
&&w(\sigma_{j}({\cal O}_{h^{N-2-j}}){\cal O}_{h^{-1-(k-N)d}}|({\cal O}_{h})^{1+(k-N)d})_{0,2|1+(k-N)d}:=\no\\
&&\frac{1}{(2 \pi \sqrt{-1})^{d+1}} \oint_{C_{0}} \frac{dz_{0}}{(z_{0})^{N}} \oint_{C_{1}} \frac{dz_{1}}{(z_{1})^{N}} \dots \oint_{C_{d}} \frac{dz_{d}}{(z_{d})^{N}}(z_{0})^{N-2-j} (z_{1} - z_{0})^{j} \biggl(\prod_{l=1}^{d} e^{k}(z_{l-1},z_{l}) \biggr)\no\\
&&\times\biggl(\prod_{l=1}^{d-1} \frac{1}{kz_{l} (2z_{l} - z_{l-1} - z_{l+1})}\biggr)\frac{1}{(z_{d})^{1+(k-N)d}}\left( d + \frac{z_0}{z_1 - z_0} \right)^{1+(k-N)d} \hspace{1cm}(2\leq N\leq k).
\label{residue2}
\ea
In the above formulas, $e^{k}(z,w)$ is given by $\prod_{j=0}^{k}(jz+(k-j)w)$, and the operation $\frac{1}{2 \pi \sqrt{-1}} \oint_{C_{i}} dz_{i}$ means taking residues at $z_{i} = 0$ for $i = 0,d$ and at $z_{i} = 0, \frac{z_{i-1} + z_{i+1}}{2}$ for $i = 1 , \ldots , d-1$. Residue integral is taken in ascending order with respect to the subscript of $z_{i}$'s. 
\label{wdefi}
\end{defi} 
In the above definition, we assume that the integer $j$ can take any non-negative integers. 

The first one, $w(\sigma_{(N-k)d+j-1}({\cal O}_{h^{N-2-j}}){\cal O}_{h^{0}})_{0,2}$, 
is given as an intersection number of the moduli space of quasimaps $\widetilde{Mp}_{0,2}(N,d)$ from $CP^{1}$ with two marked points $0,\infty\in CP^{1}$ to $CP^{N-1}$ \cite{Jin1, Jin3, S}, 
if $0\leq j\leq N-2$. \footnote{The symbol $\sigma_{j}$ means the $j$-th power of Mumford Morita class defined as the first Chern class of the line bundle on $\widetilde{Mp}_{0,2}(N,d)$
whose fiber is given as the cotangent space of $CP^{1}$ at the first marked point $0\in CP^{1}$.} 
In this case, we can express the intersection number by using elements of Chow ring of $\widetilde{Mp}_{0,2}(N,d)$:
\ba
&&w(\sigma_{(N-k)d+j-1}({\cal O}_{h^{N-2-j}}){\cal O}_{h^{0}})_{0,2}=\no\\
&&\int_{\widetilde{Mp}_{0,2}(N,d)}(H_{1}-H_{0})^{(N-k)d+j-1}(H_{0})^{N-2-j}\biggl(\prod_{i=1}^{d}\frac{e^{k}(H_{i-1},H_{i})}{kH_{i}}\biggr)(kH_{d}).
\ea
In the above formula, we interpret $\frac{e^{k}(H_{j-1},H_{j})}{kH_{j}}$ as $\prod_{i=1}^{k}(iH_{j-1}+(k-i)H_{j})$ and $H_{0},H_{1},\ldots, H_{d}$ are generators of Chow ring of $\widetilde{Mp}_{0,2}(N,d)$ that satisfy the following relations \cite{S}:
\ba
(H_{0})^{N}=0,\;(H_{j})^{N}(2H_{j}-H_{j-1}-H_{j+1})=0\;\;(j=1,2,\ldots,d-1),\;(H_{d})^{N}=0.
\ea
The factor $\prod_{l=1}^{d-1} \frac{1}{(2z_{l} - z_{l-1} - z_{l+1})}$ in (\ref{residue1}) and (\ref{residue2}) comes from the second relation $(H_{j})^{N}(2H_{j}-H_{j-1}-H_{j+1})=0$.
If $j>N-2$, we can no longer express  $w(\sigma_{(N-k)d+j-1}({\cal O}_{h^{N-2-j}}){\cal O}_{h^{0}})_{0,2}$ in terms of Chow ring because negative power of $H_{0}$ appears. But the residue integral representation (\ref{residue1}) may give us non-vanishing rational number even in this case. 

The second one,  $w(\sigma_{j}({\cal O}_{h^{N-2-j}}){\cal O}_{h^{-1-(k-N)d}}|({\cal O}_{h})^{1+(k-N)d})_{0,2|1+(k-N)d}$ is more exotic. The symbol ``$h$'' originally means hyperplane class 
in $H^{1,1}(CP^{N-1},\mathbb{C})$, but in notation of the intersecion number, negative power of $h$ appears. It is formally interpreted as a $2+(1+(k-N)d)$ pointed intersection number
of the moduli space of quasimaps  $\widetilde{Mp}_{0,2|(1+(k-N)d)}(N,d)$ constructed in \cite{JS}.
By allowing negative power of $h$ formally, this  intersection number can alternatively be represented as follows:
\ba
&&w(\sigma_{j}({\cal O}_{h^{N-2-j}}){\cal O}_{h^{-1-(k-N)d}}|({\cal O}_{h})^{1+(k-N)d})_{0,2|1+(k-N)d}=\no\\
&&\sum_{i=0}^{\mbox{min.}\{1+(k-N)d, j\}}{1+(k-N)d\choose i}d^{1+(k-N)d-i}w(\sigma_{j-i}(h^{N-2-j+i}){\cal O}_{h^{-1-(k-N)d}})_{0,d}.
\ea 
In the above formula, we assumed Hori's equation \cite{hori} for $2+m$ pointed intersection numbers:
\begin{equation}
w(\sigma_j ({\cal O}_{h^a} ) {\cal O}_{h^b} | ({\cal O}_h )^m )_{0,2|m} \stackrel{\mathrm{formally}}{=} d \cdot w(\sigma_j ({\cal O}_{h^a} ) {\cal O}_{h^b} | ({\cal O}_h )^{m-1} )_{0,2|m-1} + w(\sigma_{j-1} ({\cal O}_{h^{a+1}} ) {\cal O}_{h^b} | ({\cal O}_h )^{m-1} )_{0,2|m-1},
\end{equation}
and applied it iteratively. This equation is proved in the case of $m=1$ in \cite{JM2}. By allowing the following ``formal'' expression:
\ba
&&w(\sigma_{j-i}({\cal O}_{h^{N-2-j+i}}){\cal O}_{h^{-1-(k-N)d}})_{0,d}\stackrel{\mathrm{formally}}{=}\no\\
&&\int_{\widetilde{Mp}_{0,2}(N,d)}(H_{1}-H_{0})^{j-i}(H_{0})^{N-2-j+i}\biggl(\prod_{i=1}^{d}\frac{e^{k}(H_{i-1},H_{i})}{kH_{i}}\biggr)(kH_{d})\frac{1}{(H_{d})^{1+(k-N)d}},
\ea 
we reach the formula (\ref{residue2}). $w(\sigma_{j}({\cal O}_{h^{N-2-j}}){\cal O}_{h^{-1-(k-N)d}}|({\cal O}_{h})^{1+(k-N)d})_{0,2|1+(k-N)d}$ may also turn out to be non-vanishing for 
any non-negative integer $j$.

 In this paper, we prove the following two theorems on these numbers.
\begin{theorem}
If $N>k\geq 1$, the following equality holds.
\ba
\frac{1}{k}w(\sigma_{(N-k)d+j-1}({\cal O}_{h^{N-2-j}}){\cal O}_{h^{0}})_{0,2}=\frac{1}{j!}\frac{\d^{j}}{\d \varepsilon^{j}}\left(\frac{\prod_{r=1}^{kd}(r+k\varepsilon)}{\prod_{r=1}^{d}(r+\varepsilon)^N}\middle) \right|_{\varepsilon=0}.
\label{mainth1}
\ea
\label{main1}
\end{theorem}
\begin{theorem} If $2\leq N\leq k $, the following equality holds.
\ba
\frac{1}{k}w(\sigma_{j}({\cal O}_{h^{N-2-j}}){\cal O}_{h^{-1-(k-N)d}}|({\cal O}_{h})^{1+(k-N)d})_{0,2|1+(k-N)d}
=\frac{1}{j!}\frac{\d^{j}}{\d \varepsilon^{j}}\left(\frac{\prod_{r=1}^{kd}(r+k\varepsilon)}{\prod_{r=1}^{d}(r+\varepsilon)^N}\middle) \right|_{\varepsilon=0}.
\label{mainth2}
\ea
\label{main2}
\end{theorem}
These two theorems are extentions of our former result given in \cite{JM1}, which realized generalized hypergeometric series used in mirror computation of genus $0$ Gromov-
Witten invariants of Calabi-Yau hypersurface in $CP^{N-1}$ as a generating function of the intersection number $w(\sigma_{j}({\cal O}_{h^{N-2-j}}){\cal O}_{h^{-1}})_{0,d}$ of $\widetilde{Mp}_{0,2}(N,d)$, 
to the case of degree $k$ hypersurface in $CP^{N-1}$. Theorem \ref{main1} corresponds to Fano $(k< N)$ case, and Theorem \ref{main2} corresponds to Calabi-Yau and general 
type $(k\geq N)$ cases. 

In Fano case, Givental considered the following differential equation:
\ba
\left(\biggl(\frac{d}{dx}\biggr)^{N-1}-ke^{x}\prod_{j=1}^{k-1}(k\frac{d}{dx}+j)\right)w(x)=0.
\label{givd}
\ea
Linear independent solutions of the above equation are given as follows.
\ba
w_{j}(x)=\sum_{d=0}^{\infty}\frac{\d^{j}}{\d \varepsilon^{j}}\left(\frac{\prod_{r=1}^{kd}(r+k\varepsilon)}{\prod_{r=1}^{d}(r+\varepsilon)^N}e^{d+\varepsilon}x\middle) \right|_{\varepsilon=0}\;\;(j=0,1,\ldots,N-2).
\label{wgiv}
\ea  
In \cite{giv}, Givental computed gravitational Gromov-Witten invariant $\langle\sigma_{(N-k)d+j-1}({\cal O}_{h^{N-2-j}}){\cal O}_{h^{0}}\rangle_{0,d}$, which is defined as intersection number of 
moduli space of stable maps $\overline{M}_{0,2}(CP^{N-1},d)$, by using localization technique invented by Kontsevich \cite{km}, and proved the following theorem:
\begin{theorem}{\bf (Givental, Theorem 9.1 in \cite{giv})}\footnote{To be precise, the theorem given here is arranged by the authors from Givental's original statement.}
If $N-k\geq 2\;(k\geq 1)$, the following equality holds.
\ba
\frac{1}{k}\langle\sigma_{(N-k)d+j-1}({\cal O}_{h^{N-2-j}}){\cal O}_{h^{0}}\rangle_{0,2}=\frac{1}{j!}\frac{\d^{j}}{\d \varepsilon^{j}}\left(\frac{\prod_{r=1}^{kd}(r+k\varepsilon)}{\prod_{r=1}^{d}(r+\varepsilon)^N}\middle) \right|_{\varepsilon=0}.
\ea
\label{givth1}
\end{theorem}
Therefore, Theorem \ref{main1} corresponds to quasimap version of Theorem \ref{givth1}. Since we are treating the moduli space of quasimaps $\widetilde{Mp}_{0,2}(N,d)$, the 
equality (\ref{mainth1}) holds in the $N-k=1$ case. Origin of this difference is expalined in \cite{Jin2}.  In contrast to complexity of the proof of Theorem \ref{givth1}, due to complicated combinatorial structure of boundaries of the moduli space of stable maps, our proof of Theorem \ref{main1} is quite straghtforward and simple. 
 
 In the general type case, we can still consider the differential equation (\ref{givd}) and the series given in (\ref{wgiv}) are still formal solutions. But as was suggested in \cite{iritani},
convergence radii of these series are equal to $0$. Therefore, Theorem \ref{main2} should be regarded as a ``formal'' result. Exotic characteristics of the intersection number $w(\sigma_{j}({\cal O}_{h^{N-2-j}}){\cal O}_{h^{-1-(k-N)d}}|({\cal O}_{h})^{1+(k-N)d})_{0,2|1+(k-N)d}$ may come from this formality. Theorem \ref{main2} can be interpreted as a kind of completion of the equality observed in \cite{Jin2}:
\ba
\frac{1}{k}\cdot d^{1+(k-N)d}\cdot w({\cal O}_{h^{N-2}}{\cal O}_{h^{-1-(k-N)d}})_{0,2}=\frac{(kd)!}{(d!)^{N}}.
\ea 
 In closing this section, we mention new feature of the proof of the main thoerems, presented in Subsection 2.1. This technique drastically simplifies computational processes of 
 the proof. Hence the proof given in Subsections 2.2 and 2.3 can be regarded as simplification of the proof given in our former literature \cite{JM1}.

\vspace{2em}
{\bf Acknowledgment} 
We would like to thank Prof.~G.~Ishikawa and Prof.~A.~Tsuchida for kind encouragement. 
Our research is partially supported by JSPS grant No. 22K03289.  

\section{Proof of the Main Theorems }

\subsection{The ``Infinitesimal Displacement'' of a Pole}

In this subsection, in order to compute the residue integrals (\ref{residue1}) and (\ref{residue2}) effectively, we introduce technique of reduction of order of a pole in the residue integrals. Let $\alpha$ be any complex constant. Let $f(z, w)$ be a complex function of two variables that has the form:
\begin{equation}
f(z,w) = \frac{g(z,w)}{2w - z - \alpha}. \label{foundform}
\end{equation}
In (\ref{foundform}), $g(z,w)$ is a holomorphic function on the open subset 
\begin{equation}
B_{r_1 , r_2} := \left\{ (z,w) \in \mathbb{C}^{2} \ ; \ |z| < 2r_1 , |2w - z - \alpha| < 2r_2 \right\}
\end{equation}
for some positive real constants $r_1 , r_2$ satisfying $0 < r_1 < 2r_2 - r_1$. Moreover, let $C(0)$ and $C(\frac{z + \alpha}{2})$ be contours $z(t) := r_1 \exp (2 \pi \sqrt{-1} t) \ (r_1 > 0 \; ; \; 0 \le t \le 1)$ on $z$-plane and $w(t) := \frac{z + \alpha}{2} + r_{2} \exp (2 \pi \sqrt{-1} t) \ (r_2 > 0 \; ; \; 0 \le t \le 1)$ on $w$-plane, respectively. 
 
We consider the following residue integral:
\begin{equation}
I_j := \frac{1}{(2\pi \sqrt{-1})^2} \oint_{C_{0}} \frac{dz}{z} \oint_{C_{\frac{z + \alpha}{2}}} dw f(z, w) \left( \frac{w-z}{z} \right)^j \quad (j = 0, 1, \dots ), \label{originint}
\end{equation}
where $\frac{1}{2 \pi \sqrt{-1}} \oint_{C_{0}} dz$ and $\frac{1}{2 \pi \sqrt{-1}} \oint_{C_{\frac{z + \alpha}{2}}} dw$ are the operations of taking residue at $z = 0$ and $w = \frac{z + \alpha}{2}$, respectively. We remark here that these are realized as contour integrals $\frac{1}{2 \pi \sqrt{-1}} \oint_{C(0)} dz$ and $\frac{1}{2 \pi \sqrt{-1}} \oint_{C(\frac{z + \alpha}{2})} dw$. In (\ref{originint}), residue integrals are done from left to right in accordance with the notation used in Definition \ref{wdefi}. Hence we integrate the $z$-variable first.  The integrand in (\ref{originint}) have a higher order pole at $z=0$.  In such case, we have to compute higher derivatives with respect to the variable $z$. 
In order to avoid computing higher derivatives, we introduce the generating function of $I_{j}$'s (this operation leads to ``infinitesimal displacement'' of the pole at $z=0$). Then we can reduce our computation to taking residue of a simple pole of the $z$-variable.  Let $F(\varepsilon)$ be the generating function of $I_j \ (j=0, 1, \dots)$ given as follows:
\begin{align}
F(\varepsilon) &:= \sum_{j=0}^{\infty} I_j \varepsilon^j \notag \\
&=\frac{1}{(2\pi \sqrt{-1})^2} \sum_{j=0}^{\infty} \oint_{C_0} \frac{dz}{z} \oint_{C_{\frac{z + \alpha}{2}}} dw \;f(z, w) \left( \frac{w-z}{z} \varepsilon \right)^j \notag \\
&=\frac{1}{(2\pi \sqrt{-1})^2} \sum_{j=0}^{\infty} \oint_{C(0)} \frac{dz}{z} \oint_{C(\frac{z + \alpha}{2})} dw \;f(z, w) \left( \frac{w-z}{z} \varepsilon \right)^j , \label{genfunc}
\end{align}
where $\varepsilon$ is a small parameter. The part of $z$-integration of the above generating funcion: 
\begin{align}
G_j (w; \varepsilon) &:= \frac{1}{2 \pi \sqrt{-1}} \oint_{C(0)} \frac{dz}{z} \; f(z,w) \left( \frac{w-z}{z} \varepsilon \right)^j \notag \\
&= \frac{1}{2 \pi \sqrt{-1}} \oint_{C(0)} \frac{dz}{z} \; \frac{g(z,w)}{2w - z - \alpha} \left( \frac{w-z}{z} \varepsilon \right)^j \notag \\
&= \frac{1}{2 \pi \sqrt{-1}} \int_{0}^{1} \frac{g \left( r_1 e^{2\pi \sqrt{-1} t} ,w \right)}{2w - r_1 e^{2\pi \sqrt{-1}t} - \alpha} \left( \frac{w - r_1 e^{2\pi \sqrt{-1}t}}{r_1 e^{2\pi \sqrt{-1}t}} \varepsilon \right)^j \cdot \frac{2 \pi \sqrt{-1} r_1 e^{2 \pi \sqrt{-1}t}dt}{r_1 e^{2\pi \sqrt{-1}t}} \notag \\
&= \int_{0}^{1} \frac{g \left( r_1 e^{2\pi \sqrt{-1} t} ,w \right)}{2w - r_1 e^{2\pi \sqrt{-1}t} - \alpha} \left( \frac{w - r_1 e^{2\pi \sqrt{-1}t}}{r_1 e^{2\pi \sqrt{-1}t}} \varepsilon \right)^j dt,
\end{align}
is holomorphic for $w$ on $B_w := \left\{ w \in \mathbb{C} \ ; \ r_1 < |2w - \alpha| < 2r_2 - r_1 \right\}$.\footnote{If $w \in B_w$, then $|2w - r_1 e^{2 \pi \sqrt{-1} t} - \alpha| \le |2w - \alpha | + r_1 < (2r_2 - r_1 ) + r_1 = 2r_2$ (i.e., $(r_1 e^{2 \pi \sqrt{-1}t}, w) \in B_{r_1 , r_2}$) and $|2w - r_1 e^{2\pi \sqrt{-1} t} -\alpha | \ge |2w - \alpha | - r_1 > r_1 - r_1 = 0$.} By using Weierstrass M-test, we can easily see that we can exchange order of integration and summation in (\ref{genfunc}):
\begin{align}
F(\varepsilon) &= \frac{1}{(2\pi \sqrt{-1})^2} \oint_{C(0)} \frac{dz}{z} \oint_{C(\frac{z + \alpha}{2})} dw \sum_{j=0}^{\infty} \frac{g(z,w)}{2w - z - \alpha} \left( \frac{w-z}{z} \varepsilon \right)^j \notag \\
&= \frac{1}{1 + \varepsilon} \cdot \frac{1}{(2\pi \sqrt{-1})^2} \oint_{C(0)} dz \oint_{C(\frac{z + \alpha}{2})} dw \; \frac{g(z,w)}{2w - z - \alpha} \cdot \frac{1}{z -\frac{\varepsilon}{1 + \varepsilon} w} \notag \\
&= \frac{1}{1 + \varepsilon} \cdot \frac{1}{(2\pi \sqrt{-1})^2} \oint_{C(0)} dz \oint_{C_{\frac{z + \alpha}{2}}} dw \; \frac{g(z,w)}{2w - z - \alpha} \cdot \frac{1}{z -\frac{\varepsilon}{1 + \varepsilon} w} , \label{genint}
\end{align}
for all $\varepsilon$'s that satisfy
\begin{equation}
|\varepsilon| < m := \min \left\{ \left| \frac{r_1 e^{2 \pi \sqrt{-1} s}}{\left( \frac{r_1 e^{2 \pi \sqrt{-1} s} + \alpha}{2} + r_2 e^{2 \pi \sqrt{-1} t} \right) - r_1 e^{2 \pi \sqrt{-1}s}} \right| \; ; \; s, t \in [0, 1] \right\}.
\end{equation}
Note that this condition ensures convergence of the series $\sum_{j=0}^{\infty} \left( \frac{w-z}{z} \varepsilon \right)^j$ in (\ref{genint}). Since 
\begin{equation}
\lim_{\varepsilon \to 0} \left| \frac{\varepsilon}{1 + \varepsilon} w \right| = 0 
\end{equation}
holds and $w$ is a point belonging to the open subset $B_{r_1 , r_2}$, we can take some positive constant $r \ (< m)$ such that $\frac{\varepsilon}{1 + \varepsilon} w$ is contained in the interior of the contour $C(0)$ if $|\varepsilon| < r$. Moreover, the numerator $g(z,w)$ of the integrand in (\ref{genint}) is holomorphic on $B_{r_1 , r_2}$ that contains $C(0) \times C(\frac{z + \alpha}{2})$. Thus we can apply Cauchy's integral theorem to the $z$-integral in (\ref{genint}):
\begin{align}
&\frac{1}{1 + \varepsilon} \cdot \frac{1}{2\pi \sqrt{-1}} \oint_{C(0)} dz \; \frac{g(z,w)}{2w - z - \alpha} \cdot \frac{1}{z -\frac{\varepsilon}{1 + \varepsilon} w}  \notag \\
&= \frac{1}{1 + \varepsilon} \cdot \frac{g \left( \frac{\varepsilon}{1 + \varepsilon} w, w \right)}{2w - \frac{\varepsilon}{1 + \varepsilon} w - \alpha} \notag \\
&= \frac{1}{2 + \varepsilon} \cdot \frac{g \left( \frac{\varepsilon}{1 + \varepsilon} w, w \right)}{w - \frac{1 + \varepsilon}{2 + \varepsilon} \alpha}. 
\end{align}
Then we only have to take residue at $w = \frac{1 + \varepsilon}{2 + \varepsilon} \alpha$:
\begin{equation}
F(\varepsilon ) = \frac{1}{2 + \varepsilon} \cdot \frac{1}{2 \pi \sqrt{-1}} \oint_{C_{\frac{1 + \varepsilon}{2 + \varepsilon} \alpha}} dw \; \frac{g \left( \frac{\varepsilon}{1 + \varepsilon} w, w \right)}{w - \frac{1 + \varepsilon}{2 + \varepsilon} \alpha},
\end{equation}
where $\frac{1}{2 \pi \sqrt{-1}} \oint_{C_{\frac{1 + \varepsilon}{2 + \varepsilon} \alpha}} dw$ is the operator of taking residue at $w = \frac{1 + \varepsilon}{2 + \varepsilon} \alpha$\footnote{Formally, we have $w = \frac{\frac{\varepsilon}{1 + \varepsilon} w + \alpha}{2} \Longleftrightarrow w = \frac{1 + \varepsilon}{2 + \varepsilon} \alpha$.}. 

With these discussions, we have proved the following lemma:
\begin{lem}
Let $f(z,w)$ be a complex function of the form
\begin{equation}
f(z,w) = \frac{g(z,w)}{2w - z - \alpha}
\end{equation}
and assume that $g(z,w)$ is holomorphic on some open subset of $\mathbb{C}^2$ that contains $B_{r_1 , r_2}$ for some $r_1 , r_2 $. Then we can choose some constant $r (> 0)$ such that the following equality:
\begin{equation}
\frac{1}{(2\pi \sqrt{-1})^2} \sum_{j=0}^{\infty} \oint_{C_0} \frac{dz}{z} \oint_{C_{\frac{z + \alpha}{2}}} dw \; f(z, w) \left( \frac{w-z}{z} \right)^j \varepsilon^j = \frac{1}{2 + \varepsilon} \cdot \frac{1}{2 \pi \sqrt{-1}} \oint_{C_{\frac{1 + \varepsilon}{2 + \varepsilon} \alpha }} dw \frac{g \left( \frac{\varepsilon}{1 + \varepsilon} w ,w \right)}{w - \frac{1 + \varepsilon}{2 + \varepsilon} \alpha}, \label{lem1}
\end{equation}
holds for all $\varepsilon$'s that satisfy $|\varepsilon| < r$. In particular, the generating function $F(\varepsilon)$ of the integral $I_j$ is holomorphic at $\varepsilon = 0$.
\end{lem}

\subsection{Proof of Theorem \ref{main1}}

In this section, we prove Theorem 1 by using Lemma 1. By Definition \ref{wdefi}, $w(\sigma_{(N-k)d+j-1}({\cal O}_{h^{N-2-j}}){\cal O}_{h^{0}})_{0,2}$ is given by
\begin{align}
w(\sigma_{(N-k)d+j-1}({\cal O}_{h^{N-2-j}}){\cal O}_{h^{0}})_{0,2} &= \frac{1}{(2 \pi \sqrt{-1})^{d+1}} \oint_{C_0} \frac{dz_0}{(z_0 )^N} \dots \oint_{C_d} \frac{dz_d}{(z_d )^N} \notag \\
&\quad \cdot (z_0 )^{N-2-j} (z_1 - z_0 )^{(N-k)d+j-1} \frac{\prod_{i=1}^{d} e^k (z_{i-1} , z_{i})}{\prod_{i=1}^{d-1} kz_i (2z_i - z_{i+1} - z_{i-1})}  \notag \\
&= \frac{1}{(2 \pi \sqrt{-1})^{d+1}} \oint_{C_0} \frac{dz_0}{z_0} \oint_{C_1} \frac{dz_1}{(z_1 )^{N}} \dots \oint_{C_d} \frac{dz_d}{(z_d )^N} \notag \\
&\quad \cdot \frac{(z_1 - z_0 )^{(N-k)d-1}}{2z_1 - z_2 - z_0} \cdot \frac{e^{k}(z_0 , z_1 )}{z_0} \cdot \frac{\prod_{i=2}^{d} e^k (z_{i-1} , z_{i})}{\prod_{i=1}^{d-2} (2z_{i+1} - z_{i+2} - z_{i})} \notag \\
&\quad \cdot \frac{1}{\prod_{i=1}^{d-1} kz_i} \cdot \left( \frac{z_1 - z_0}{z_0} \right)^{j} ,
\end{align}
where
\begin{equation}
e^{k}(z,w):=\prod_{j=0}^{k}((k-j)z+jw)\;\;\;(N-k\geq 1)
\end{equation}
is a degree ($k+1$) polynomial and $\frac{1}{2 \pi \sqrt{-1}} \oint_{C_i} dz_i \ (i = 0 , \dots , d)$ is the operation of taking residue(s) at
\begin{equation}
\begin{cases}
z_i = 0 & (i= 0 , d), \\
z_i = 0 , \frac{z_{i-1} + z_{i+1}}{2} & (i = 1, \dots , d-1).
\end{cases}
\end{equation}
Note that $e^k (z , w)$ is divisible by $z$ and $w$ (and therefore $e^k (z, 0) \equiv 0$). In order to prove our assertion, we introduce the generating function of the above integrals:
\begin{align}
F_0 (\varepsilon ) &:= \sum_{j=0}^{\infty} w(\sigma_{(N-k)d+j-1}({\cal O}_{h^{N-2-j}}){\cal O}_{h^{0}})_{0,2} \; \varepsilon^j \notag \\
&= \sum_{j=0}^{\infty} \frac{1}{(2 \pi \sqrt{-1})^{d+1}} \oint_{C_0} \frac{dz_0}{z_0} \oint_{C_1} dz_1 \oint_{C_2} \frac{dz_2}{(z_2 )^{N}} \dots \oint_{C_d} \frac{dz_d}{(z_d )^N} \frac{f_0 (z_0 , \dots, z_d )}{2z_1 - z_0 - z_2} \left( \frac{z_1 - z_0}{z_0} \right)^{j} \varepsilon^j ,
\end{align}
where $f_0 (z_0 , \dots , z_d )$ is defined by 
\begin{align}
f_0 (z_0 , \dots , z_d ) &:= \frac{(z_1 - z_0 )^{(N-k)d-1}}{(z_1 )^{N}} \cdot \frac{e^k (z_0 , z_1 )}{z_0} \notag \\
&\quad \cdot \frac{\prod_{i=2}^{d} e^k (z_{i-1} , z_{i})}{\prod_{i=1}^{d-2} (2z_{i+1} - z_{i+2} - z_{i})} \cdot \frac{1}{\prod_{i=1}^{d-1} kz_i}.
\end{align}
With this set-up, we have only to prove the following equality:
\begin{equation}
F_0 (\varepsilon ) = k \cdot \frac{\prod_{r=1}^{kd}(r+k\varepsilon)}{\prod_{r=1}^{d}(r+\varepsilon)^N} \quad (\mbox{for any sufficiently small} \ \varepsilon).
\end{equation}
Note that since $e^{k} (z_{i-1} , z_{i})$ is divisible by $z_{i-1}$, $f_0 (z_0 , \dots , z_d )$ is holomorphic at the point $(z_0 , \dots , z_d )$ such that
\begin{equation}
2z_{i} - z_{i-1} - z_{i+1} \neq 0 \quad (i = 2, \dots , d-1), \quad z_1 \neq 0.
\end{equation}
Thus we can apply Lemma 1 for $\frac{f_0 (z_0 , \dots , z_d )}{2z_1 - z_0 - z_2 }$ by taking some constant $r_{0} (> 0)$. Then we obtain
\begin{align}
&\sum_{j=0}^{\infty} \frac{1}{(2 \pi \sqrt{-1})^{2}} \oint_{C_0} \frac{dz_0}{z_0} \oint_{C_1} dz_1 \frac{f_0 (z_0 , \dots , z_d )}{2z_1 - z_0 - z_2} \left( \frac{z_1 - z_0}{z_0} \right)^{j} \varepsilon^j \notag \\
&= \sum_{j=0}^{\infty} \frac{1}{(2 \pi \sqrt{-1})^{2}} \oint_{C_0} \frac{dz_0}{z_0} \oint_{C_{1, \frac{z_0 + z_2}{2}}} dz_1 \frac{f_0 (z_0 , \dots , z_d )}{2z_1 - z_0 - z_2} \left( \frac{z_1 - z_0}{z_0} \right)^{j} \varepsilon^j \notag \\
&\quad + \sum_{j=0}^{\infty} \frac{1}{(2 \pi \sqrt{-1})^{2}} \oint_{C_0} \frac{dz_0}{z_0} \oint_{C_{1,0}} dz_1 \frac{f_0 (z_0 , \dots , z_d )}{2z_1 - z_0 - z_2} \left( \frac{z_1 - z_0}{z_0} \right)^{j} \varepsilon^j \notag \\
&= \frac{1}{2 + \varepsilon} \cdot \frac{1}{2 \pi \sqrt{-1}} \oint_{C_{1,\frac{1 + \varepsilon}{2 + \varepsilon} z_2}} dz_1 \frac{f_0 \left( \frac{\varepsilon}{1 + \varepsilon} z_1 , z_1 , z_2 , \dots , z_d \right)}{z_{1} - \frac{1 + \varepsilon}{2 + \varepsilon} z_2} \notag \\
&\quad + \sum_{j=0}^{\infty} \frac{1}{(2 \pi \sqrt{-1})^{2}} \oint_{C_0} \frac{dz_0}{z_0} \oint_{C_{1,0}} dz_1 \frac{f_0 (z_0 , \dots , z_d )}{2z_1 - z_0 - z_2} \left( \frac{z_1 - z_0}{z_0} \right)^{j} \varepsilon^j \quad (|\varepsilon| < r_0 ), \label{int0}
\end{align}
where $\frac{1}{2\pi \sqrt{-1}} \oint_{C_{1,\alpha}} dz_1$ ($\alpha \in \mathbb{C}$) is the operation of taking residue at $z_1 = \alpha$.
For later use, we also deonote by $\frac{1}{2\pi \sqrt{-1}} \oint_{C_{j,\alpha}} dz_{j}$ ($\alpha \in \mathbb{C}$) the operation of taking residue at $z_j = \alpha$.
Since 
\begin{equation}
\frac{e^k \left( \frac{\varepsilon}{1 + \varepsilon} z_1 , z_1 \right)}{\frac{\varepsilon}{1 + \varepsilon} z_1} = k \left( \prod_{r=1}^k (r + k \varepsilon ) \right) \left( \frac{z_1}{1 + \varepsilon} \right)^{k} \not\equiv 0
\end{equation}
and
\begin{equation}
\left( z_1 - \frac{\varepsilon}{1 + \varepsilon} z_1 \right)^{(N-k)d-1} = \left( \frac{1}{1 + \varepsilon} z_1 \right)^{(N-k)d-1} \not\equiv 0,
\end{equation}
the 1st term of (\ref{int0}) is
\begin{align}
&\frac{1}{2 + \varepsilon} \cdot \frac{1}{2 \pi \sqrt{-1}} \oint_{C_{1, \frac{1 + \varepsilon}{2 + \varepsilon} z_2}} dz_1 \; \frac{f_0 \left( \frac{\varepsilon}{1 + \varepsilon} z_1 , z_1 , z_2 , \dots , z_d \right)}{z_1 - \frac{1 + \varepsilon}{2 + \varepsilon} z_2} \notag \\
&= \frac{k}{(1 + \varepsilon )^{(N-k)(d-1)-1} (2 + \varepsilon )} \frac{\prod_{r=1}^{k} (r + k \varepsilon )}{(1 + \varepsilon )^N} \notag \\
&\quad \cdot \frac{1}{2 \pi \sqrt{-1}} \oint_{C_{1, \frac{1 + \varepsilon}{2 + \varepsilon} z_2 }} dz_1 \; \frac{f_1 (z_1 , z_2 , \dots , z_d )}{z_1 - \frac{1 + \varepsilon}{2 + \varepsilon} z_2},
\end{align}
where
\begin{align}
f_1 (z_1 , \dots , z_d ) &:= (z_1 )^{(N-k)(d-1)-1} \cdot \frac{e^{k}(z_1 ,z_2 )}{kz_1} \cdot \frac{1}{2z_2 - z_1 - z_3} \notag \\
&\quad \cdot \frac{\prod_{i=3}^{d} e^k (z_{i-1} , z_{i})}{\prod_{i=2}^{d-2} (2z_{i+1} - z_{i+2} - z_{i})} \cdot \frac{1}{\prod_{i=2}^{d-1} kz_i}
\end{align}
and it is holomorphic at the point $(z_1 , \dots , z_d )$ such that
\begin{equation}
2z_{i+1} - z_i - z_{i+2} \neq 0 \quad (i = 1 , \dots , d-2). 
\end{equation}
On the other hand, we can compute the 2nd term of (\ref{int0}) in the same way as in the discussion in Subsection 2.1:
\begin{align}
&\sum_{j=0}^{\infty} \frac{1}{(2 \pi \sqrt{-1})^{2}} \oint_{C_0} \frac{dz_0}{z_0} \oint_{C_{1,0}} dz_1 \frac{f_0 (z_0 , \dots , z_d )}{2z_1 - z_0 - z_2} \left( \frac{z_1 - z_0}{z_0} \right)^{j} \varepsilon^j \notag \\
&= \frac{1}{2 + \varepsilon} \cdot \frac{1}{2 \pi \sqrt{-1}} \oint_{C_{1,0}} dz_1 \frac{f_0 \left( \frac{\varepsilon}{1 + \varepsilon} z_1 , z_1 , z_2 , \dots , z_d \right)}{z_{1} - \frac{1 + \varepsilon}{2 + \varepsilon} z_2} \notag \\
&= \frac{k}{(1 + \varepsilon )^{(N-k)(d-1)-1} (2 + \varepsilon )} \frac{\prod_{r=1}^{k} (r + k \varepsilon )}{(1 + \varepsilon )^N} \notag \\
&\quad \cdot \frac{1}{2 \pi \sqrt{-1}} \oint_{C_{1,0}} dz_1 \; \frac{f_1 (z_1 , z_2 , \dots , z_d )}{z_1 - \frac{1 + \varepsilon}{2 + \varepsilon} z_2} \notag \\
&=0.
\label{vani1}
\end{align}
Here, we take the integral contour of  $\oint_{C_{1,0}}dz_{1}$ as $z_{1}(t)=r_{1}\exp(2\pi\sqrt{-1}t)\;\;(r_{1}>0, 0\leq t\leq 1)$ and assume that $z_{2}$ satisfies the condition: $r_{1}< |\frac{1 + \varepsilon}{2 + \varepsilon}| |z_2|$. \footnote{Later, we impose analogous conditions on $z_{3},\cdots, z_{d}$ in evaluating $\oint_{C_{j,0}}dz_{j}\;(j=3,\cdots d)$ in order to guarantee vanishing of the terms 
arising from  $\oint_{C_{j,0}}dz_{j}$.} 
Therefore we obtain
\begin{equation}
F_0 (\varepsilon ) = \frac{k}{(1 + \varepsilon )^{(N-k)(d-1)-1} (2 + \varepsilon )} \frac{\prod_{r=1}^{k} (r + k \varepsilon )}{(1 + \varepsilon )^N} F_1 (\varepsilon ) \quad (| \varepsilon | < r_{0}), 
\end{equation}
where we set
\begin{equation}
F_1 (\varepsilon) := \frac{1}{(2 \pi \sqrt{-1})^{d}} \oint_{C_{1, \frac{1 + \varepsilon}{2 + \varepsilon} z_2}} dz_1 \oint_{C_2} \frac{dz_2}{(z_2 )^N} \dots \oint_{C_d} \frac{dz_d}{(z_d )^N} \; \frac{f_1 (z_1 , \dots , z_d )}{z_1 - \frac{1 + \varepsilon}{2 + \varepsilon} z_2}. \label{thm1int1}
\end{equation}

Next, we consider the following integration of $\frac{f_1 (z_1 , \dots , z_d )}{z_1 - \frac{1 + \varepsilon}{2 + \varepsilon} z_2}$:
\begin{align}
&\frac{1}{(2 \pi \sqrt{-1})^{2}} \oint_{C_{1, \frac{1 + \varepsilon}{2 + \varepsilon} z_2}} dz_1 \oint_{C_2} \frac{dz_2}{(z_2 )^N} \frac{f_1 (z_1 , \dots , z_d )}{z_1 - \frac{1 + \varepsilon}{2 + \varepsilon} z_2} \notag \\
&= \frac{1}{(2 \pi \sqrt{-1})^{2}} \oint_{C_{1, \frac{1 + \varepsilon}{2 + \varepsilon} z_2}} dz_1 \oint_{C_{2, \frac{z_1 + z_3}{2}}} \frac{dz_2}{(z_2 )^N} \frac{f_1 (z_1 , \dots , z_d )}{z_1 - \frac{1 + \varepsilon}{2 + \varepsilon} z_2} \notag \\
&\quad + \frac{1}{(2 \pi \sqrt{-1})^{2}} \oint_{C_{1, \frac{1 + \varepsilon}{2 + \varepsilon} z_2}} dz_1 \oint_{C_{2,0}} \frac{dz_2}{(z_2 )^N} \frac{f_1 (z_1 , \dots , z_d )}{z_1 - \frac{1 + \varepsilon}{2 + \varepsilon} z_2}. \label{int2} 
\end{align}

In the same way as the discussion in Subsection 2.1, the 1st term of (\ref{int2}) is computed as follows: 
\begin{align}
&\frac{1}{(2 \pi \sqrt{-1})^{2}} \oint_{C_{1, \frac{1 + \varepsilon}{2 + \varepsilon} z_2}} dz_1 \oint_{C_{2, \frac{z_1 + z_3}{2}}} \frac{dz_2}{(z_2 )^N} \frac{f_1 (z_1 , \dots , z_d )}{z_1 - \frac{1 + \varepsilon}{2 + \varepsilon} z_2} \notag \\
&= \frac{1}{2 \pi \sqrt{-1}} \oint_{C_{2, \frac{2 + \varepsilon}{3 + \varepsilon} z_3}} \frac{dz_2}{(z_2 )^{N}} \; f_1 \left( \frac{1 + \varepsilon}{2 + \varepsilon} z_2 , z_2 , z_3 , \dots , z_d  \right) \notag \\
&= \frac{(1 + \varepsilon )^{(N-k)(d-1)-1}}{(2 + \varepsilon )^{(N-k)(d-2)-2}(3 + \varepsilon )} \frac{\prod_{r=k+1}^{2k} (r + k \varepsilon )}{(2 + \varepsilon )^{N}} \notag \\
&\quad \cdot \frac{1}{2 \pi \sqrt{-1}} \oint_{C_{2, \frac{2 + \varepsilon}{3 + \varepsilon} z_3 }} dz_2 \; \frac{f_2 (z_2 , \dots , z_d )}{z_2 - \frac{2 + \varepsilon}{3 + \varepsilon} z_3},
\end{align}
where we defined
\begin{equation}
f_2 (z_2 , \dots , z_d ) := (z_2 )^{(N-k)(d-2)-1} \cdot \frac{\prod_{i=3}^{d} e^k (z_{i-1} , z_{i})}{\prod_{i=2}^{d-2} (2z_{i+1} - z_{i+2} - z_{i})} \cdot \frac{1}{\prod_{i=2}^{d-1} kz_i}.
\end{equation}
 Then the function $f_2 (z_2 , \dots , z_d )$ is holomorphic at the point $(z_2 , \dots , z_d )$ where the following conditions are satisfied:
\begin{equation}
2z_{i+1} - z_i - z_{i+2} \neq 0 \quad (i = 2 , \dots , d-2). 
\end{equation}
On the other hand, the 2nd term of (\ref{int2}) vanishes in the same way as the computation in (\ref{vani1}):
\begin{align}
&\frac{1}{(2 \pi \sqrt{-1})^{2}} \oint_{C_{1, \frac{1 + \varepsilon}{2 + \varepsilon} z_2}} dz_1 \oint_{C_{2,0}} \frac{dz_2}{(z_2 )^N} \frac{f_1 (z_1 , \dots , z_d )}{z_1 - \frac{1 + \varepsilon}{2 + \varepsilon} z_2} \notag \\
&= \frac{(1 + \varepsilon )^{(N-k)(d-1)-1}}{(2 + \varepsilon )^{(N-k)(d-2)-2}(3 + \varepsilon )} \frac{\prod_{r=k+1}^{2k} (r + k \varepsilon )}{(2 + \varepsilon )^{N}} \notag \\
&\quad \cdot \frac{1}{2 \pi \sqrt{-1}} \oint_{C_{2,0}} dz_2 \; \frac{f_2 (z_2 , \dots , z_d )}{z_2 - \frac{2 + \varepsilon}{3 + \varepsilon} z_3} \notag \\
&= 0.
\end{align}
Here, we take the integral contour of  $\oint_{C_{2,0}}dz_{2}$ as $z_{2}(t)=r_{2}\exp(2\pi\sqrt{-1}t)\;\;(r_{2}>0, 0\leq t\leq 1)$ and assume that $r_{2}$ and $z_{3}$ satisfiy
the conditions: $r_{1}< |\frac{1 + \varepsilon}{2 + \varepsilon}| r_{2}$, ${2}<|\frac{2 + \varepsilon}{3 + \varepsilon}| |z_3|$, respectively.  Hence we have
\begin{align}
F_1 (\varepsilon ) &= \frac{(1 + \varepsilon )^{(N-k)(d-1)-1}}{(2 + \varepsilon )^{(N-k)(d-2)-2}(3 + \varepsilon )} \frac{\prod_{r=k+1}^{2k} (r + k \varepsilon )}{(2 + \varepsilon )^{N}} F_2 (\varepsilon ) \quad (| \varepsilon | < r_{1}), \label{int1}
\end{align}
where
\begin{equation}
F_2 (\varepsilon ) :=  \frac{1}{(2 \pi \sqrt{-1})^{d-1}} \oint_{C_{2, \frac{2 + \varepsilon}{3 + \varepsilon} z_3}} dz_2 \oint_{C_3} \frac{dz_3}{(z_3 )^N} \dots \oint_{C_d} \frac{dz_d}{(z_d )^N} \; \frac{f_2 (z_2 , \dots , z_d )}{z_2 - \frac{2 + \varepsilon}{3 + \varepsilon} z_3}.
\end{equation}
By repeating the procedures so far, we reach the following
expression:
\begin{align}
F_0 (\varepsilon ) &= \sum_{j=0}^{\infty} w(\sigma_{(N-k)d+j-1}({\cal O}_{h^{N-2-j}}){\cal O}_{h^{0}})_{0,2} \; \varepsilon^j \notag \\
&= \frac{k}{(1 + \varepsilon )^{(N-k)(d-1)-1} (2 + \varepsilon )} \frac{\prod_{r=1}^{k} (r + k \varepsilon )}{(1 + \varepsilon )^N} F_1 (\varepsilon ) \notag \\
&= \frac{k}{(2 + \varepsilon )^{(N-k)(d-2)-1} (3 + \varepsilon )} \frac{\prod_{r=1}^{2k} (r + k \varepsilon)}{\prod_{r=1}^{2} (r + \varepsilon )^{N}} F_2 (\varepsilon ) \notag \\
&= \dotsb \notag \\
&= \frac{k}{(d-1 + \varepsilon )^{(N-k) \cdot 1 - 1} (d + \varepsilon )} \frac{\prod_{r=1}^{(d-1)k} (r + k \varepsilon )}{\prod_{r=1}^{d-1} (r + \varepsilon )^{N}} F_{d-1} (\varepsilon) \quad (\mbox{for any sufficiently small} \ \varepsilon),
\end{align}
where
\begin{equation}
F_{d-1} (\varepsilon) := \frac{1}{(2 \pi \sqrt{-1})^2} \oint_{C_{d-1, \frac{d-1 + \varepsilon}{d + \varepsilon} z_{d}}} dz_{d-1} \oint_{C_d} \frac{dz_d}{(z_d )^{N}} \frac{(z_{d-1})^{(N-k) \cdot 1 - 1}}{z_{d-1} - \frac{d - 1 + \varepsilon}{d + \varepsilon} z_{d}} \frac{e^{k}(z_{d-1} , z_{d})}{kz_{d-1}}.
\end{equation}
Then we can easily evaluate this integral as
\begin{equation}
F_{d-1} (\varepsilon ) = \frac{(d - 1 + \varepsilon )^{(N-k) \cdot 1 - 1}}{(d + \varepsilon )^{-1}} \frac{\prod_{r = (d-1)k + 1}^{kd} (r + k \varepsilon )}{(d + \varepsilon )^{N}}.
\end{equation}
In this way, we finally obtain
\begin{equation}
F_0 (\varepsilon ) = k \cdot \frac{\prod_{r=1}^{kd}(r+k\varepsilon)}{\prod_{r=1}^{d}(r+\varepsilon)^N} \quad (\mbox{for any sufficiently small} \ \varepsilon),
\end{equation}
which completes the proof of Theorem \ref{main1}. \ $\square$

\subsection{Proof of Theorem \ref{main2}}
As was done in the proof of Theorem 1, we consider the generating function:
\begin{align}
G_0 (\varepsilon) &:= \sum_{j=0}^{\infty} w(\sigma_{j}({\cal O}_{h^{N-2-j}}){\cal O}_{h^{-1-(k-N)d}}|({\cal O}_{h})^{1+(k-N)d})_{0,2|1+(k-N)d} \; \varepsilon^j.
\end{align}
By using (\ref{residue2}) in Definition \ref{wdefi}, $G_{0}(\varepsilon)$ is given as the following residue integral:
\begin{align}
G_0 (\varepsilon)&= \frac{1}{(2 \pi \sqrt{-1})^{d+1}} \sum_{j=0}^{\infty} \oint_{C_0} \frac{dz_0}{z_0} \oint_{C_1} dz_1 \oint_{C_2} \frac{dz_2}{(z_2 )^N} \dots \oint_{C_d} \frac{dz_d}{(z_d )^N} \; \frac{g_0 (z_0 , \dots , z_d )}{2z_1 - z_0 - z_2} \left( \frac{z_1 - z_0}{z_0} \right)^{j} \varepsilon^{j},
\end{align}
where we set $g_0 (z_0 , \dots , z_d )$ as
\begin{align}
g_0 (z_0 , \dots , z_d ) &:= \frac{1}{(z_1 )^{N}} \cdot \left( d + \frac{z_0}{z_1 - z_0} \right)^{1+(k-N)d} \cdot \frac{e^k (z_0 , z_1 )}{z_0} \notag \\
&\quad \cdot \frac{\prod_{i=2}^{d} e^k (z_{i-1} , z_{i})}{\prod_{i=1}^{d-2} (2z_{i+1} - z_{i+2} - z_{i})} \cdot \frac{1}{\prod_{i=1}^{d-1} kz_i} \cdot  (z_d )^{-1-(k-N)d}.
\end{align}
Since $g_0 (z_0 , \dots , z_d  )$ is holomorphic at the point $(z_0 , \dots , z_d )$ such that
\begin{equation}
2z_{i} - z_{i-1} - z_{i+1} \neq 0 \quad (i = 2, \dots , d-1), \quad z_1 \neq 0 , \quad z_d \neq 0,
\end{equation}
we can apply Lemma 1 and the remaining processes go in the same way as the proof of Theorem 1. $\Box$

\end{document}